\documentclass[10pt]{article}
\usepackage{amssymb}
\usepackage{amsmath}
\usepackage{amsfonts}
\usepackage{amsthm}

\usepackage[english]{babel}
\usepackage[affil-it]{authblk}

\newtheorem{theorem}{Theorem}

\newtheorem{proposition}{Proposition}
\newtheorem{corollary}{Corollary}

\newtheorem{conjecture}{Conjecture}
\title{MMS-type problems for Johnson scheme\thanks{The reported study was funded by RFBR according to the research project 18-31-00126 }}

\author{I.Yu.Mogilnykh, K.V.Vorob'ev, A.A.Valyuzhenich, %
  \thanks{E-mail address: \texttt{ivmog@math.nsc.ru, vorobev@math.nsc.ru, graphkiper@mail.ru}}}
\affil{Sobolev Institute of Mathematics, Novosibirsk State University,
 Novosibirsk, Russia}
\begin{document}

\maketitle

\begin{abstract}
 In the current work we consider the
minimization problems for the number of nonzero or negative values
of vectors from the first and second eigenspaces of the Johnson
scheme respectively. The topic is a meeting point for
generalizations of the Manikam-Mikl\'{o}s-Singhi conjecture proven in
\cite{Blin} and the minimum support problem for the eigenspaces of
the Johnson graph, asymptotically solved in \cite{MVV}.
\end{abstract}

\section{Introduction}

 Let $V_i$ be $i$th eigenspace of a
symmetric association scheme $(X,\{R_0,\ldots,R_d\})$, $i\in
\{0,\ldots,d\}$. Following \cite{Bier1}, given an eigenvector
$v\in V_i$ denote by $X_{+}(v)=\{x\in X:v_x>0\}$, $X_{-}(v)=\{x\in
X:v_x<0\}$, $X_{0}(v)=\{x\in X:v_x=0\}$. A pair of $w$-subsets are
in $i$th relation, if their intersection is of size $w-i$. The
$w$-subsets of $\{1,\ldots, n\}$ together with $w+1$ relations
above define {\it the Johnson scheme} and the first relation
defines the Johnson graph $J(n,w)$.
 The eigenvalues of the Johnson
scheme are known as the values of the Eberlein polynomials
$E_k(i,w,n)=\sum_{j=0}^{k}(-1)^j{i \choose j}{w-i \choose k-j}{n-w-i \choose k-j}$,
$k,i\in \{0,1,\ldots w\}$. For this scheme (graph) by $V_i$, we
denote the eigenspace corresponding to the eigenvalue
$\lambda_i(n,w)=E_1(i,w,n)=(w-i)(n-w-i)-i$ for $i\in \{0,1,\ldots
w\}$.

 In the
current correspondence we consider the following two
characteristics for the Johnson scheme $J(n,w)$:
$$m_i^-(n,w)=\underset{v:v\in V_i, X_0(v)= \emptyset}{min}|X_{-}(v)|,$$
$$m_i^{0}(n,w)=\underset{v:v\in V_i, v \not\equiv 0}{min}|X_{+}(v)|+|X_{-}(v)|.$$

When $i=1$, the first number was suggested to be ${n-1 \choose w-1}$,
for $n\geq4w$, which is known as the Manikam-Micklos-Singhi
conjecture \cite{Mani}, \cite{Mani2}. The problem was completely
solved in the paper of Blinovski \cite{Blin} recently, after works
with quadratic \cite{Alon} and linear improvements \cite{Pok}.

For a vector $v$ by the support of this vector we mean the value $X_{+}(v)+X_{-}(v)$. The number $m_i^{0}(n,w)$ was shown to be equal to $2^i{n-2i
\choose w-i}$ (along with the description of vectors attaining the
bound) for sufficiently large $n$ in \cite{MVV}. In the paper we
focus on the case when $i=1$ and show that for any $n$ and $w$ the
minimum of the support of vectors from the first eigenspace of
$J(n,w)$ is attained on the vectors from two classes having rather
simple structure (see Section 3).

 Bier and Delsarte \cite{BierDelsarte} proposed to investigate the invariant
 $\underset{v:v\in V_i, X_0(v)= \emptyset}{min}|X_{-}(v)|$ for
 classical association schemes with further generalizations where $v$ is from the
 direct sum of several eigenspaces. They obtained several bounds involving such well-known combinatorial
 concepts
as coverings, completely regular codes, additive codes and
designs. The current study is motivated by a recent progress in
the area of completely regular codes and equitable partitions in
Johnson graphs. In particular, the characterizations of equitable
2-partitions of $J(n,3)$ in \cite{GG} for odd $n$, completely regular codes in $J(n,w)$ with eigenvalue $\lambda_2$
having nontrivial distance were characterized by Martin in
\cite{Martin}.

An eigenvector $u$ of the antipodal Johnson graph $J(2w,w)$
corresponding to $\lambda_i$ is such that its absolute values on
the pairs of antipodal pairs are equal and signs are the same or
opposite depending on the parity of $i$ \cite{BCN}[p. 142-143].
So, in case of odd $i$ we have that $m_i^{-}(2w,w)={2w-1 \choose
w-1} $.

In \cite{BierDelsarte} it was shown that
\begin{equation}\label{BDbounds}\frac{(^{n}_{w})}{|D|}\leq m_i^{-}(n,w)\leq |C|,\end{equation} where $C$ and
$D$ are codes (subsets of the vertices of $J(n,w)$), whose
characteristic functions belong to $V_i\oplus V_0$ and
$\underset{j\in \{0,\ldots,w\}\setminus i}{\oplus}V_j$
respectively. Eliminating a constant vector from the
characteristic function of $C$, we see that there is a two-valued
eigenvector $v$ from $V_i$ such that $v_x\neq v_y$ iff $x\in C$,
$y\notin C$. In other words, $(C,\overline{C})$ is an equitable
2-partition of $J(n,w)$. Suppose that there is a
$(w-1)-(n,w,1)$-design $C$ then its size is the value for
$m_{w}^{-}(n,w)$. Indeed, such a design produces an equitable
2-partition $(C, \overline{C})$ of $J(n,w)$ with eigenvalue
$\lambda_w(n,w)$, see \cite{Martin}. On the other hand the
"anticode" $D$ could be chosen to be the set $\{x:y\subset x\}$
where $y$ is a $w-1$-element subset. The set $D$ is a Delsarte
clique in the Johnson graph which is a completely regular code
with eigenvalues $\lambda_0,\ldots,\lambda_{w-1}$, so the
characteristic function of $D$ is orthogonal to $V_w$
\cite{Bier2}. The smallest open case is $i=2,w=3$, since
$m_1^{-}(n,2)$ was shown to be $\lceil n/2\rceil$ in \cite{Bier2}.
Again, for $n=1,3 (mod\mbox{ }6)$, $w=3$ the best known
 "anticode" $D$ from $V_0\oplus V_1 \oplus V_3$ is
a Steiner triple system, so from (\ref{BDbounds}) we have that

\begin{equation}\label{LBn3}n-2\leq m_2^{-}(n,3).\end{equation}

The bound (\ref{LBn3}) could be tightened up to $2n-9$ by
considering a modification of the weight distribution lower bound
\cite{KMP} with a generalization for arbitrary $w$, which we
discuss in Section 4.1. The choice of $C$ in (\ref{BDbounds}) is
generalized to be a part of an equitable partition with appropriate eigenvalue. This gives an upper bound in case of
$J(n,3)$ for odd $n$ and $i=2$ (see Section 4.2), where no equitable 2-partitions
exist \cite{GG}. For even $n$ the upper bound
(\ref{BDbounds}) from equitable 2-partitions of $J(n,3)$ is
$n(n-2)/2$.

\section{Definitions and Preliminaries}


\subsection{Equitable partitions}

Let $G$ be an undirected graph. An {\it equitable $r$-partition}
with parts $C_1,\ldots, C_{r}$ of the vertex set of $G$ is called
equitable if for any $i,j \in \{1,\ldots,r\}$ a vertex from $C_i$
has exactly $A_{ij}$ neighbors in $C_j$. The matrix
$A=(A_{ij})_{i,j\in \{1,\ldots,r\}}$ is called the {\it quotient
matrix}. An eigenvalue of the quotient matrix $A$ is called {\it
an eigenvalue} of the partition. Given an eigenvector $u$ of $A$
corresponding to an eigenvalue $\lambda$ define $u^G$ to be the
vector, indexed by the vertices of $G$ such that $u^G_x=u_i$, if
$x\in C_i$. The vector $u^G$ is an eigenvector of the adjacency
matrix of $G$ corresponding to $\lambda$ \cite{CDZ}[\S 4.5].
 In view of the said above, the upper bound in (\ref{BDbounds}) is generalized as follows:

 \begin{proposition}\label{Ub}
Let $u$ be an eigenvector without zero entries of the quotient
matrix of an equitable partition of the Johnson graph $J(n,w)$
with parts $C_1,\ldots,C_r$. Then
$$m_i^{-}(n,w)\leq \sum_{j: u_j<0}|C_j|.$$
 \end{proposition}





\subsection{The first eigenspace of $J(n,w)$}

Consider the eigenvectors of the complete graph $K_n=J(n,1)$
indexed by integers from $\{1,\ldots,n\}$. The graph has two
eigenvalues: $n-1$ and $-1$. An eigenvector
$a=(\alpha_1,\ldots,\alpha_n)$ of the graph corresponding to the
eigenvalue $-1$ could be characterized as a solution for the
equation: $\alpha_1+\ldots+\alpha_n=0$. The first eigenspace of
the Johnson graph $J(n,w)$ could be obtained from that of the
complete graph using the isomorphism established by the inclusion
mapping $I$ (see \cite{Delsarte}): the image $I_{}(a)$ is such
that $(I(a))_x=\underset{i\in x}{\sum} \alpha_i$.


Consider the following
two equitable 2-partitions of $J(n,w)$: $(\{x:1 \in x\}$, $\{x:1
\notin x\})$ and $(\{x:2 \in x\}$, $\{x:2 \notin x\})$ \cite{Martin}.
Denote by $v^{1,2}$ the difference of two eigenvectors of $J(n,w)$
arising from these partitions, i.e. $v^{1,2}_x=0,$ if $1,2$ are
simultaneously in or are not in $x$, $v^{1,2}_x=1$, if $1\in x$,
$2\notin x$, $v^{1,2}_x=-1$, if $1\notin x$, $2\in x$. In \cite{MVV}
it was shown that the minimum support eigenvectors from the first
eigenspace are exactly $v^{1,2}$ up to appropriate permutation of coordinate positions starting with large enough $n$ (as
well as a generalization of the result for any eigenspace). It is
easy to see that $v^{1,2}$ is $I(e_1-e_2)$, where $e_1, e_2$ are
$1$-st and $2$-nd vectors of the standard basis.

In Section 3 we extend results from \cite{MVV} in further details.
We show that for any $n$ the minimum support eigenvector is either
$v^{1,2}$ or $I(a)$, where $a$ is a two-valued $(-1)$-eigenvector
of $J(n,1)$.

\section{Minimum support $\lambda_1$-eigenvectors}

\begin{theorem}{\label{th1}}
Let $v$ be $\lambda_1$-eigenvector of $J(n,w)$, $n\geq 2w$, $w\ge 2$ with minimum support. Then $v$ is $I(e_1-e_2)$ or $I(\sum_{i=1}^{k}{e_i}-\frac{k}{n-k}\sum_{i=k+1}^{n}{e_i})$ for some $k\in \{2,3,\dots,n-2\}$ such that
 $\frac{kw}{n}\in \mathbb{N}$ up to a permutation of coordinate positions and the multiplication by a scalar.
In particular,

$$m_1^0(n,w)=\min{(2{n-2 \choose w-1},{n \choose w}-\max_{k\in \{2,3,\dots,n-2\},\frac{kw}{n}\in \mathbb{N}}{ {k \choose \frac{kw}{n}}{n-k \choose \frac{(n-k)w}{n}}  })} $$

\end{theorem}
\proofname

As it was mentioned above every $\lambda_1$-eigenvector equals $I(a)$ for some vector $a=(\alpha_1,\ldots,\alpha_n)$
such that $\alpha_1+\ldots+\alpha_n=0$. Our next goal is to determine values $\alpha_1,\ldots,\alpha_n$ for which the support of $I(a)$ is minimal. Let $v=I(a)$ be a $\lambda_1$-eigenvector with minimum support. Since the vector $I(e_1-e_2)$ has the size of the support equal $2{n-2 \choose w-1}$, we shall assume that the size of the support of $I(a)$ is not more than $2{n-2 \choose w-1}$. Let us denote by $m$ the size of $\{\alpha_1,\ldots,\alpha_n\}$.
There are two different cases:

\begin{enumerate}
\item[$m\ge 3.$]
Without loss of generality we can assume that $\alpha_1,\alpha_2$
and $\alpha_3$ are pairwise different. Take arbitrary subsets
$A_1,A_2$ of the set  $\{4,\ldots,n\}$ of cardinalities $w-1$ and
$w-2$ respectively. Clearly, there are at least $2$ nonzero values
among $I(a)_{i\cup A_1}=\alpha_i+\sum_{k\in A_1}\alpha_k, i=1,2,3$
and at least $2$ nonzero values among $I(a)_{\{i,j\}\cup
A_2}=\alpha_i+\alpha_j+\sum_{k\in A_2}\alpha_k,$ $i,j \in
\{1,2,3\}$, $i\ne j$. So the support of $v$ is at least $2{n-3
\choose w-1}+2{n-3 \choose w-2}=2 {n-2 \choose w-1}$. By
hypothesis the vector $I(a)$ has minimal size of the support, so
we conclude that $I(a)_x=0$ for any $w$-subset $x$ of
$\{4,\ldots,n\}$. In other words, $I'(a')$ is the zero vector,
where $a'$ is obtained from $a$ by removing its first 3 entries,
$I'$ is the inclusion mapping from $J(n-3,1)$ to $J(n-3,w)$.

We have that $\sum_{i=4,\ldots,n}\alpha_i=\frac{\sum_{x\subset
\{4,\ldots,n\}, |x|=w}I'(a')_x}{(^{n-4}_{w-1})}=0$. Therefore the
vector $a'=(\alpha_4,\ldots,\alpha_n)$ belongs to $V_1(n-3,1)$ and
is the zero vector because $I'(a')$ is the zero vector and $I'$ is
an isomorphism from $V_1(n-3,1)$ to $V_1(n-3,w)$.

   From the
above there are exactly $2$ nonzero values among
$\alpha_1+\alpha_2, \alpha_1+\alpha_3, \alpha_2+\alpha_3$.
Consequently, we can consider $\alpha_3=0$ and
$\alpha_1=-\alpha_2$, which means that $v$ is equal to
$cI(e_1-e_2)$ for some constant $c$.
\item [$m=2.$]
 Without loss of generality we can take $\alpha_1=\alpha_2=\ldots=\alpha_k=\hat{\alpha}$ and $\alpha_{k+1}=\alpha_{k+2}=\ldots =\alpha_{n}=\hat{\beta}$ for some integer $k$ such that $2 \leq k\ \leq n-2$. Using the equality $\alpha_1+\ldots+\alpha_n=0$ we have $\hat{\beta}=-\hat{\alpha}\frac{k}{n-k}$.
 Let us take an arbitrary vertex $x$ of $J(n,w)$. It is easy to see that $I(a)_x=0$ if and only if $x$ has exactly $\frac{kw}{n}$ ones in first $k$ coordinate positions and $w-\frac{kw}{n}$ ones in the rest $n-k$ coordinate positions. Particularly, $\frac{kw}{n}$ must be an integer.
 Therefore, the support of $v$ equals ${n \choose w}- {k \choose \frac{kw}{n}}{n-k \choose \frac{(n-k)w}{n}}$. Taking the minimum of this expression over all admissible $k$ we obtain the statement of the theorem.

\end{enumerate}
In the case $m=1$ we automatically obtain the all-zero eigenvector $v$, which is not possible.

\qed

The theorem $1$ reduces the problem of minimizing $m_1^0(n,w)$ to
the comparison of two expressions containing binomial
coefficients.

\section{Bounds on $m_i^{-}(n,w)$}

\subsection{A lower bound on $m_i^{-}(n,w)$}

Let $v$ be an eigenvector of the Johnson graph $J(n,w)$ without
zero entries, $x$ be the vertex such that $v_x$ is negative and
takes maximum absolute value over all negative entries of $v$.
Consider the distance partition $(C_0,\ldots,C_w)$ of the vertices
of $J(n,w)$ with respect to the vertex $x$. It is well-known that
the sum of the entries of $v$ on $C_k$ is expressed using the
Eberlein polynomials and the value $v_x$:
$$\sum_{y \in C_k}v_y=v_xE_k(i,w,n).$$ Let $E_k(i,w,n)$ be non-negative. Then by the choice of $v_x$ with the maximum absolute value we see that there are at
least $|E_k(i,w,n)|$ negative values for $v_y$ in $C_k$. Moreover,
 there are more than $|E_k(i,w,n)|$ negative $v_y$'s not
 less then $v_x$, because there is at least one positive $v_y$
 in $C_k$, since obviously $|E_k(i,w,n)|<|C_k|$ for $k>0$. Thus we obtain the following bound.

 \begin{theorem}
$m_i^{-}(n,w)\geq 1+\underset{k>0:
E_k(i,w,n)\ge 0}{\sum}(|E_k(i,w,n)|+1)$.
 \end{theorem}
The consideration for the proof above is similar to the one for
the weight distribution bound on the number of nonzeros for the
eigenvector of distance-regular graph, see \cite{KMP}. The values
of the Eberlein polynomials for $i=2$ and $w=3$ are as follows
$E_0(2)=1,$ $E_1(2)=n-7$, $E_2(2)=11-2n$, $E_3(2)=n-5$. Therefore,
we have the bound below.

\begin{corollary}{\label{lowbound}}
$m_2^{-}(n,3)\geq 2n-9$.
\end{corollary}

\subsection{An upper bound on $m_2^{-}(n,3)$}
Let $n$ be $2r$. The following construction could be found in
\cite{GP} (see also \cite{AM}). Consider the complement of a
perfect matching on vertices labeled with $\{1,\ldots,2k\}$ to a
complete bipartite graph. Then the triples of vertices are parted
into three orbits $C_1$, $C_2$, $C_3$ with respect to the action
of the automorphism group of the graph. The triples of $C_1$
consist of vertices belonging to the same part, the triples of
$C_2$ induce a walk of length 2 in the graph, the triples of $C_3$
contain exactly one pair of adjacent vertices. Any two parts could
be merged and result in equitable 2-partition of triples, e.g. the
Johnson graph $J(2r,3)$ \cite{AM}. In particular, the partition
($C_1'=C_1\cup C_2$, $C_2'=C_3$) has the following quotient matrix:
$$\left(%
\begin{array}{cc}
  3(2r-5) & 6 \\
  4(r-2) & 2r-1 \\
\end{array}%
\right)$$
 whose eigenvalues are $\lambda_0(n,3)$ and
$\lambda_2(n,3)$. The parts are in $4(r-2)$ to 6 ratio, so in view
of the Proposition \ref{Ub} we see that $$m_2^{-}(n,3)\leq
n(n-2)/2.$$

Let $n$ be $2r+1$. Consider the graph $G$ with $2r+1$ vertices
which is a union of an isolated vertex and the graph $G'$
consisting the complement of a perfect matching on vertices of
size $r$ to a complete bipartite graph. We have the following
orbits of triples of vertices:

\begin{description}
    \item[$C_1$] The vertices of the triple are in one part of $G'$
    \item[$C_2$] The vertices of the triple induce a walk of length 2 in $G'$
    \item[$C_3$] The vertices of the triple belong to $G'$ and
    contain
    only two adjacent vertices
    \item[$C_4$] Two nonadjacent vertices belong to different parts of $G'$ and the
    third one is isolated
    \item[$C_5$] Two vertices are in one part of $G'$ and the third one
    is isolated
    \item[$C_6$] Two vertices are adjacent and the third one is isolated
\end{description}

The equitable partition $(C_1,\ldots, C_6)$ of $J(2r+1,3)$ has the
following quotient matrix:

$$\left(%
\begin{array}{cccccc}
  3(r-3) & 3(r-2) & 6 & 0 & 3 & 0 \\
  r-2 & 5r-13 & 6 & 0 & 1 & 2 \\
  r-2 & 3(r-2) & 2r-1 & 1 & 1 & 1 \\
  0 & 0 & 2(r-1) & 0 & 2(r-1) & 2(r-1) \\
  r-2 & r-2 & 2 & 2 & 2(r-2) & 2(r-1) \\
  0 & 2(r-2) & 2 & 2 & 2(r-1) & 2(r-2) \\
\end{array}%
\right).$$ The matrix has eigenvector $(3,3,4-2r,2-2r,1,1)$
corresponding to eigenvalue $\lambda_2(2r+1,3)=2r-6$. By
Proposition \ref{Ub}, we see that $$m_2^{-}(n,3)\leq
|C_3|+|C_4|=2r(r-1)+r=(n-1)(n-2)/2.$$

Thus we obtain
\begin{theorem}{\label{upbound}}

$$m_2^{-}(n,3)\leq
\begin{cases}
n(n-2)/2,&\text{if $n$ is even;}\\
(n-1)(n-2)/2,&\text{if $n$ is odd.}

\end{cases}  $$

\end{theorem}

\section{Conclusion}

Theorem 1 reduces the problem of finding $m_1^0(n,w)$ to the
determination which one of values
$$ \left( {n \choose w}-\max_{k\in \{2,3,\dots,n-2\},\frac{kw}{n}\in \mathbb{N}}{ {k \choose \frac{kw}{n}}{n-k \choose \frac{(n-k)w}{n}}  } \right)  \text{ or }   2{n-2 \choose w-1}$$
is smaller. In \cite{MVV} it was shown that the second one is the
answer starting from some value $n_0(w)$. We have compared these
values for $6 \leq n \leq 600$ and $3\leq w \leq \frac{n}{2}$ and
consequently found corresponding $m_1^0(n,w)$. Based on these
computational results we state the following conjecture:

\begin{conjecture}
For $w \geq 5$ and $n\geq 2w+1$ the following identity holds
$$m_1^0(n,w)=2{n-2 \choose w-1}.$$
\end{conjecture}

For $w<5$ we have found several curious examples:
\begin{enumerate}
\item $m_1^0(6,2)=6$ is attained on the vector $v=I(e_1+e_2+e_3-e_4-e_5-e_6)$,

\item $m_1^0(8,2)=12$ is attained on vectors $v=I(e_1+e_2+e_3+e_4-e_5-e_6-e_7-e_8)$ and $u=(e_1-e_2)$,

\item $m_1^0(9,3)=39$ is attained on the vector $v=I(2e_1+2e_2+2e_3-e_4-e_5-e_6-e_7-e_8-e_9)$,

\item $m_1^0(10,4)=110$ is attained on the vector $v=I(e_1+e_2+e_3+e_4+e_5-e_6-e_7-e_8-e_9-e_{10})$.
\end{enumerate}

Let us notice that it is not hard to show using theorem \ref{th1}
and basic properties of binomial coefficients that
$$m_1^0(2w,w)={2w \choose w}-2{2w-2 \choose w-1}$$ which is
attained on the vector $I((w-1)(e_1+e_2)-\sum_{i=3}^{2w}{e_i})$.

In the theorem \ref{upbound} we described a construction providing
a quadratic on $n$ upper bound for the characteristic
$m_2^{-}(n,3)$. At the same time, the corollary \ref{lowbound}
gives us a lower bound which is linear on $n$. The real behaviour
of the growth rate of $m_2^{-}(n,3)$ remains to be an intriguing
open problem.

The characteristic $\underset{v:v\in V_i, X_0(v)=
\emptyset}{min}|X_{-}(v)|$ considered by Bier and Delsarte
\cite{BierDelsarte} requires that $v$ does not have zero entries.
It may be interesting in the future research to remove this
condition and try to find this value in this case for classical
association schemes.


\section{Acknowledgements}
The authors thank the anonymous referee of the paper \cite{MVV}
for pointing out that there is a connection between the problem of
minimizing of Johnson scheme eigenvector's support and the
Manikam-Mikl\'{o}s-Singhi conjecture, which led to the current
research.

\end{document}